\documentclass[11pt]{amsart}
\usepackage{amsfonts}
\usepackage{amsmath,amssymb}
\usepackage[foot]{amsaddr}
\usepackage{amsthm}
\usepackage{tikz}
\usetikzlibrary{arrows.meta,positioning,calc,decorations.pathreplacing}
\usepackage{bm,bbm}
\usepackage{color}

\usepackage{bigints}

\usepackage{mathtools}
\usepackage{commath}

\numberwithin{equation}{section}

\allowdisplaybreaks

\theoremstyle{plain}
\newtheorem{theorem}{Theorem}[section]
\newtheorem{lemma}[theorem]{Lemma}
\newtheorem{proposition}[theorem]{Proposition}
\newtheorem{corollary}[theorem]{Corollary}
\theoremstyle{definition}

\newtheorem{remark}[theorem]{Remark}

\def\ms{\mathsf}
\def\beqn{\begin{equation}}
\def\beqn*{$$}
\def\eeqn{\end{equation}}

\newcommand{\bx}{{\bf x}}
\newcommand{\by}{{\bf y}}
\newcommand{\bz}{{\bf z}}

\def\P{\mathbb{P}}

\def\E{\mathbb{E}}
\def\Pn{\mathcal P_n}

\newcommand{\reals}{{\mathbb R}}

\newcommand{\R}{\reals}

\newcommand{\vep}{\varepsilon}
\newcommand{\bbz}{\protect{\mathbb Z}}

\newcommand{\X}{{\mathcal{X}}}
\newcommand{\Y}{{\mathcal{Y}}}

\newcommand{\B}{\mathcal B}

\newcommand{\ta}{\theta}

\newcommand{\one}{{\mathbf 1}}

\newcommand{\remove}[1]{}

\newcommand{\Z}{\mathcal Z}

\def\enp{\end{proof}}
\def\bel{\begin{lemma}}
\def\bep{\begin{proof}}
\def\enl{\end{lemma}}

\def\ms{\mathsf}

\newcommand{\Pnr}{\mathcal P_n|}

\newcommand{\QmM}{Q_\ell \setminus M_\ell}
\newcommand{\bbT}{{\mathbb T}}
\newcommand{\al}{\alpha}

\begin{document}

\bibliographystyle{abbrv}

\renewcommand{\baselinestretch}{1.05}

\title[LDP in the dense regime]
{Large deviation principles for component counts in the dense regime}

\author{Christian Hirsch}
\address{Department of Mathematics, Aarhus University, Ny Munkegade 118, 8000 Aarhus C, Denmark}
\email{christian.hirsch@math.au.dk}

\author{Takashi Owada}
\address{Department of Statistics, Purdue University, West Lafayette, IN 47907, USA}
\email{owada@purdue.edu}

\thanks{Owada's research was partially supported by the AFOSR grants FA9550-22-1-0238 at Purdue University. }



\begin{abstract}
This paper develops large deviation principles for component counts in dense
random geometric graphs.  We consider a homogeneous Poisson point process on
the flat \(d\)-dimensional torus and study the number of isolated components
of fixed order \(k\) at connection radii above the connectivity scale.  Two
types of large deviation behaviour are investigated.  First, we prove a
Donsker--Varadhan large deviation principle in the regime where the expected
number of such components tends to infinity.  Second, we study the opposite
extreme regime, where this expectation tends to zero, and establish
fixed order tail measure limits for the associated \(k\)-component
isolation-volume point process.  As a consequence, we obtain functional large
deviations for the moving threshold component count process.  
\end{abstract}

\maketitle

\section{Introduction}

Random geometric graphs are a fundamental model for spatial networks and
continuum connectivity, with roots in stochastic geometry and continuum
percolation; see, for instance, the monographs of Penrose~\cite{penrose:2003}
and Meester and Roy~\cite{meester:roy:1996}.  They also appear naturally in
coverage and clump count problems~\cite{hall:1986,hall:1988}, wireless-network
models~\cite{baccelli:blaszczyszyn:2009}, and geometric statistics through
nearest neighbor functionals and graph based multivariate testing procedures
\cite{friedman:rafsky:1979,bickel:breiman:1983,henze:1988}.  In this work, we investigate the large deviation behavior of the $k$-component count in the dense regime.

There are two mathematical difficulties in the analysis of component counts in random geometric graphs.  The first is that large
deviations for component structure are delicate even outside the geometric
setting.  In sparse Erd\H{o}s--R\'enyi and sparse inhomogeneous random graphs,
large deviation principles (LDPs) have been developed for empirical component
profiles, including microscopic, mesoscopic and macroscopic components and
type-refined component data
\cite{andreis:konig:patterson:2020,andreis:konig:langhammer:patterson:2023}.
For spatial random graphs, related large deviation questions have also been
studied for the largest cluster in finite boxes, including long-range
percolation, geometric inhomogeneous random graphs and age-dependent random
connection models~\cite{jorritsma:komjathy:mitsche:2024}.  These works illustrate the difficulties of establishing large deviation
principles for component structure.  They do not, however, cover the dense random geometric graph
considered here.

The second difficulty is geometric. 
Large deviations in spatial probabilistic models often lead to 
condensation phenomena.  In random geometric graphs, upper tails for the edge
count are known to be realized by concentrating excess edges in a small spatial
region~\cite{chatterjee:harel:2020}.  Related condensation phenomena also arise
in large deviation analyses of spatial many-particle systems, for example in
symmetrized Brownian bridge models connected with Bose--Einstein condensation
\cite{adams:konig:2008}.  The mechanism in the present paper is different, but
it has the same broad feature that the rare event is realized through a
structured spatial arrangement rather than through homogeneous fluctuations.

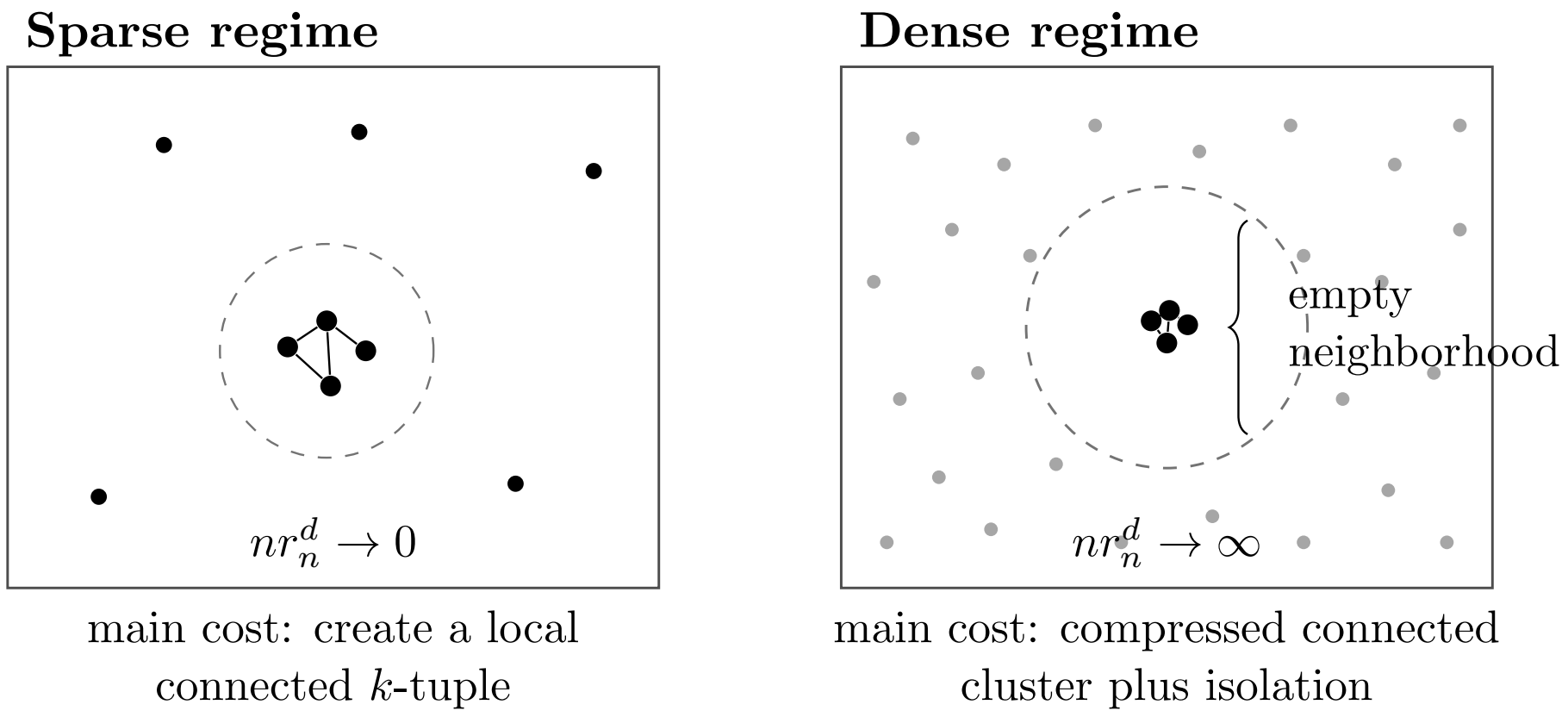
\begin{figure}[t]
\centering
\begin{tikzpicture}[
        scale=1,
        every node/.style={font=\small},
        point/.style={circle, fill=black, inner sep=1.25pt},
        faintpoint/.style={circle, fill=black!35, inner sep=1.05pt},
        clusterpoint/.style={circle, fill=black, inner sep=1.6pt},
        edge/.style={black, line width=0.45pt},
        radius/.style={black!55, dashed, line width=0.45pt},
        box/.style={draw=black!70, line width=0.55pt}
]

% Left panel: sparse regime
\begin{scope}
\node[anchor=west,font=\bfseries] at (0,4.25) {Sparse regime};

\draw[box] (0,0) rectangle (5,4);

% background sparse points
\node[point] at (0.7,0.7) {};
\node[point] at (1.2,3.4) {};
\node[point] at (3.9,0.8) {};
\node[point] at (4.5,3.2) {};
\node[point] at (2.7,3.5) {};

% local k tuple
\node[clusterpoint] (a1) at (2.15,1.85) {};
\node[clusterpoint] (a2) at (2.45,2.05) {};
\node[clusterpoint] (a3) at (2.75,1.82) {};
\node[clusterpoint] (a4) at (2.48,1.55) {};

\draw[edge] (a1) -- (a2);
\draw[edge] (a2) -- (a3);
\draw[edge] (a2) -- (a4);
\draw[edge] (a1) -- (a4);

\draw[radius] (2.45,1.82) circle (0.82);

\node[align=center] at (2.5,-0.55)
{main cost: create a local\\connected \(k\)-tuple};

\node[align=center] at (2.5,0.35)
{\(nr_n^d\to0\)};
\end{scope}

% Right panel: dense regime
\begin{scope}[xshift=6.4cm]
\node[anchor=west,font=\bfseries] at (0,4.25) {Dense regime};

\draw[box] (0,0) rectangle (5,4);

% dense background points outside exclusion zone
\foreach \x/\y in {
0.35/0.35,0.75/0.85,1.15/0.45,1.65/0.95,2.15/0.35,2.85/0.55,3.55/0.35,4.2/0.75,4.65/0.35,
0.45/1.45,1.05/1.65,3.85/1.45,4.55/1.65,
0.25/2.35,0.85/2.75,1.45/2.55,3.55/2.55,4.15/2.35,4.75/2.75,
0.55/3.45,1.25/3.25,1.95/3.55,2.75/3.35,3.45/3.55,4.25/3.25,4.75/3.55}
{
  \node[faintpoint] at (\x,\y) {};
}

% empty surrounding neighborhood
\draw[fill=white, draw=black!55, dashed, line width=0.55pt] (2.5,2.0) circle (1.08);

% compressed cluster
\node[clusterpoint] (b1) at (2.38,2.05) {};
\node[clusterpoint] (b2) at (2.52,2.13) {};
\node[clusterpoint] (b3) at (2.66,2.02) {};
\node[clusterpoint] (b4) at (2.50,1.88) {};

\draw[edge] (b1) -- (b2);
\draw[edge] (b2) -- (b3);
\draw[edge] (b2) -- (b4);
\draw[edge] (b1) -- (b4);

\draw[decorate,decoration={brace,amplitude=4pt},line width=0.45pt]
      (3.12,1.18) -- (3.12,2.82)
      node[midway,right=5pt,align=left] {empty\\neighborhood};

\node[align=center] at (2.5,-0.55)
{main cost: compressed connected\\cluster plus isolation};

\node[align=center] at (2.5,0.35)
{\(nr_n^d\to\infty\)};
\end{scope}

\end{tikzpicture}
\caption{
Two mechanisms for fixed size components in random geometric graphs.  In the
sparse regime, the main cost is creating a local connected \(k\)-tuple.  In the
dense regime, local connectivity is abundant, and the cost is instead to form a
compressed connected \(k\)-cluster whose connection neighborhood contains no
external points.  }
\label{fig:sparse-dense-mechanism}
\end{figure}

For finite components in random geometric graphs with connectivity radius $r_n$, a sparse regime
large deviation theory is available when \(nr_n^d\to0\)
\cite{hirsch:owada:2023}.  In that regime, the obstruction is to create local
connected \(k\)-tuples at all; once such a tuple appears, separation from the
rest of the process is comparatively inexpensive.  The natural speed is
therefore of sparse combinatorial order.  The dense regime considered in the
present paper has a different obstruction.  Since \(nr_n^d\to\infty\), local
connectivity is abundant.  The rare event is that a microscopic connected
cluster remains isolated from the surrounding dense cloud of points. See Figure \ref{fig:sparse-dense-mechanism} for an illustration.

This dense obstruction is related to  dense
nearest neighbor and low degree phenomena.  For \(k\)-nearest neighbor balls,
the relevant one point probability can be estimated using Poisson lower tails,
and dense large deviation principles are available in Euclidean and hyperbolic
settings~\cite{hirsch:kang:owada:2023,hirsch:otto:owada:thale:2025}.  For
\(k\)-components, the corresponding one point probability is more difficult to estimate.  The \(k\) points must
form a connected cluster, but the cluster must also be isolated from all other
points.  Even the first- and second-order asymptotics of this probability in
the dense regime require high-intensity component asymptotics
\cite{penrose:yang:2025}, based on the compression mechanism for finite
clusters in continuum percolation~\cite{alexander:1993}.

\medskip
\noindent
The high-intensity component asymptotic of Penrose--Yang is the starting
point for our LDP.  Write
\[
        a_n:=\ta_dnr_n^d\to\infty,
        \qquad
        m_k:=(k-1)(d-1),
\]
where $\ta_d$ is the volume of a unit ball in $\R^d$, 
and assume that
\(
        b_n:=n a_n^{-m_k}e^{-a_n}\to\infty. 
\)
The relevant quantity is the Palm probability \(p_k(\lambda)\) that the
origin belongs to a component of order \(k\) in the random geometric graph
\(G(\mathcal H_\lambda\cup\{0\},1)\), with unit connectivity radius, where
\(\mathcal H_\lambda\) is a homogeneous Poisson point process of intensity
\(\lambda\) on \(\R^d\):
\[
        p_k(\lambda)\sim
        \al_k\lambda^{-m_k}e^{-\ta_d\lambda},
        \qquad \lambda\to\infty .
\]
The polynomial correction and the constant \(\al_k\) reflect the compression
mechanism for finite clusters.  Our first main result turns this local
high-intensity component estimate into an LDP for the
global count.  Theorem~\ref{t:LDP.Skn} proves an LDP for the component counts of order $k$ with speed \(b_n\) and rate function
\[
        I_k(x):=x\log(x/s_k)-x+s_k,\qquad x\ge0,
\]
for some constant $s_k>0$. 
Both the speed and the constant in the rate function are therefore inherited
from the dense cluster formation probability.

\medskip
\noindent
The same scale separates two complementary asymptotic regimes.  If
$b_n\to\infty$, 
then components of order \(k\) occur often enough for their total count to
exhibit a many-component large deviation behavior.  This is the regime of
Theorem~\ref{t:LDP.Skn}.  On the other hand, if $b_n\to0$,
then components of order \(k\) are rare at the relevant scale.  In
this regime, the natural object is no longer a normalized count but the point
process of centered external isolation radii of \(k\)-components.  Our second
main result gives the corresponding fixed order tail measure limits in the
sense of \(\mathbb M\)-convergence.

The form of the rate function and its proof show that, on the LDP scale, the dense \(k\)-components behave as
if they were generated by an approximately Poisson process of rare local
clusters.  We partition the torus at spatial scale \(b_n^{-1/d}\), introduce a
trimmed blocked count \(\eta_{k,n}\), and prove that the one-block count
converges to a Poisson random variable with mean \(s_k\).  A diluted-cube
argument gives the uniform exponential integrability needed for convergence of
logarithmic moment generating functions.  The LDP for the blocked count then
follows from the G\"artner--Ellis theorem, with the lower bound at zero checked
directly.  Finally, an exponential equivalence argument shows that components
intersecting block boundaries are negligible on the scale \(b_n\).  The
blocking and exponential equivalence strategy follows the dense spatial LDP
approach of~\cite{hirsch:kang:owada:2023}, while the one box
probability estimate follows from the high-intensity component asymptotic
of~\cite{penrose:yang:2025}.

The fixed radius LDP studies one deterministic slice of the radius-indexed
component process. 
In many situations one would like
a functional version, where several thresholds can be studied simultaneously.
For nearest neighbor distances, this type of rare event asymptotic can be
formulated through convergence of exceedance point processes
\cite{hirsch:kang:owada:2023}.  In the second part of the paper, we develop
the component analogue in the tail measure theory of Dombry, Tillier and
Wintenberger~\cite{dombry:tillier:wintenberger:2022}.  This approach records
the first order exceedance measure and distinguishes rare event strata according
to the number of extreme components.  Configurations with one, two, or more
surviving microscopic components occur on different large deviation scales.

The contributions of the paper are as follows.
\begin{enumerate}
\item We prove a dense LDP for fixed order component
counts in random geometric graphs.  In the regime \(a_n=\ta_dnr_n^d\to\infty\)
and \(b_n=n a_n^{-m_k}e^{-a_n}\to\infty\), the normalized $k$-component count
 satisfies an LDP with speed \(b_n\) and Poissonian rate
function \(I_k(x)=x\log(x/s_k)-x+s_k\). 

\item We identify the high-intensity \(k\)-component asymptotic of
Penrose--Yang as determining the large deviation scale.  Thus,
the speed and the rate function constant are determined by the probability
that \(k\) points form a compressed connected component and remain isolated
from the surrounding dense cloud.

\item We prove a rare event tail measure theorem for the associated
\(k\)-component distance.  This places the component process within the
tail measure theory of Dombry, Tillier and Wintenberger.  The higher order
convergence distinguishes configurations according to the number of
extreme points: the limit concerns configurations with at least \(q\)
surviving components and yields the corresponding \(q\)-th order tail measure.
Thus, the result describes higher order rare events in addition to the
first order Poisson exceedance limit.
\end{enumerate}

The rest of the paper is organized as follows.  Section~\ref{sec:LDP}
introduces the model, establishes the high-intensity component constants,
and proves the LDP through blocking, one-block Poisson approximation,
uniform exponential integrability and exponential equivalence.
Section~\ref{sec:extreme-k-clusters} proves
the tail measure theorem for the \(k\)-component distance.

\section*{Acknowledgements}

The authors thank T.~Kang for fruitful discussions during the
early stages of this project.

The mathematical results and original proofs in
Section~\ref{sec:LDP} were developed independently by the
authors, without generative AI assistance. During subsequent
revisions, ChatGPT (OpenAI) was used merely to improve its exposition. 

Generative AI was used more extensively in
Section~\ref{sec:extreme-k-clusters}. In particular, ChatGPT
assisted in exploring alternative proof strategies, including
a direct approach through factorial moments and a simplification
of the interaction estimates using translation invariance.
It also suggested combining  several auxiliary arguments
and assisted with the exposition.

All AI-assisted mathematical suggestions incorporated into the
manuscript were independently checked and verified by the
authors. The authors take full responsibility for the mathematical
statements, proofs, and conclusions.

\section{Large deviation principle}  \label{sec:LDP}

We establish the LDP in  the dense regime, i.e., $nr_n^d\to \infty$ as $n\to\infty$. Fix $k,d\ge2$ and define 
$$
h(x_1,\dots, x_k) := \one \big\{   G \big( \{x_1,\dots, x_k\}, 1 \big) \text{ is connected} \big\},  \ \ \ x_i \in \R^d, 
$$
where $G \big( \{x_1,\dots, x_k\}, 1 \big) $ denotes a geometric graph on a vertex set $\{ x_1,\dots, x_k \}$ with unit connectivity radius; that is, an edge is present between $x_i$ and $x_j$ whenever the distance between $x_i$ and $x_j$ is less than one. %$\|x_i-x_j\|\le 1$. 
The scaled version of $h$ is defined by 
\begin{equation}  \label{e:def.hn}
h_n(x_1,\dots,x_k) := h(r_n^{-1}x_1,\dots,r_n^{-1}x_k),  \ \ \ x_i\in\R^d. 
\end{equation}
Denote by $\Pn$  the homogeneous Poisson point process with intensity $nf$ where $f$ is constant uniform on the flat torus $\bbT^d:=  [0,1]^d/\sim$, which  is equipped with the toroidal metric
$\ms{dist} (x,y) = \min_{z\in \bbz^d} \| x-y+z \|$.  
For a finite point set $\Y$ in $\R^d$ or $\bbT^d$ and $r>0$, write
\begin{equation}\label{e:def.union.ball}
\B(\Y;r):=\bigcup_{p\in\Y}B_r(p),
\end{equation}
where $B_r(p)$ is the closed ball of radius $r$ in the underlying
Euclidean or toroidal metric.  For a $k$-point set $\Y\subset\R^d$
and a finite external configuration $\xi\subset\R^d\setminus\Y$, define
\[
c(\Y,\xi):=\one\{\xi(\B(\Y;1))=0\}.
\]
Define also 
$$
c_n(\Y, \xi) := c(r_n^{-1}\Y, r_n^{-1}\xi), \ \ \text{ and } \ \ g_n(\Y, \xi) := h_n(\Y)c_n(\Y, \xi). 
$$
%For $r>0$, define
%\[
%\begin{aligned}
%        h_r(x_1,\dots,x_k)&:=h(r^{-1}x_1,\dots,r^{-1}x_k),
%        \qquad c_r(\Y,\xi):=c(r^{-1}\Y,r^{-1}\xi),\\
%        g_r(\Y,\xi)&:=h_r(\Y)c_r(\Y,\xi).
%\end{aligned}
%\]
%Let $S_{k,n}(r)$ denote the number of connected components of order $k$ in
%$G(\Pn,r)$; equivalently,
%\[
%        S_{k,n}(r):=
%        \sum_{\Y\subset \Pn,\, |\Y|=k} g_r(\Y,\Pn\setminus \Y).
%\]
Let $S_{k,n}$ denote the number of connected components of order $k$ in $G(\Pn, r_n)$; that is, 
\begin{equation}  \label{e:def.Tkn}
S_{k,n}
:= \sum_{\Y\subset \Pn, \, |\Y|=k} g_n(\Y, \Pn\setminus \Y), 
\end{equation}
where $(r_n)_{n\ge1}$ satisfies  $nr_n^d\to\infty$ and tends to $0$ faster than the threshold for vanishing of  components of order $k$. Our setup and proof techniques are closely related to \cite{penrose:yang:2025, iyer:yogeshwaran:2020}. Recall that the threshold radius $r_n$  in \cite{penrose:yang:2025} is of the form 
\begin{equation}  \label{e:critical.radius}
r_n=\bigg( \frac{\log n -m_k\log \log n +C}{n\ta_d} \bigg)^{1/d}, 
\end{equation}
for some constant $C\in \R$, where $m_k:=(k-1)(d-1)$ and  $\ta_d$ denotes the volume of a unit ball in $\R^d$. 
In particular, under the radius \eqref{e:critical.radius}, the count of $k$-components in a random geometric graph converges weakly to a Poisson random variable~\cite[Theorem 5.1]{penrose:yang:2025}. 
Since our interest lies in the LDP in the dense regime, we consider the radius 
\begin{equation}  \label{e:LDP.radius}
r_n = \bigg( \frac{a_n}{n\ta_d} \bigg)^{1/d}.
\end{equation}
We assume that $(a_n)_{n\ge1}$ is positive, tends to infinity, and satisfies
\begin{equation}  \label{e:growth.an}
a_n-\log n +m_k\log \log n \to-\infty,  \ \ \ n\to\infty. 
\end{equation}
The speed sequence for our LDP is defined as
$        b_n:=n a_n^{-m_k}e^{-a_n}$ for $n\ge1$.
Note that $b_n\to\infty$ as $n\to\infty$. For the proof, write 
\begin{equation}  \label{e:rewrite.bn}
b_n = e^{\log n -m_k\log \log n -a_n} \Big( \frac{\log n}{a_n} \Big)^{m_k}.
\end{equation}
By \eqref{e:growth.an}, we have $a_n^{-1}\log n \ge 1/2$ for all large
$n$, and thus, $b_n\to\infty$ follows. 

Finally, we  introduce
an ``integer" block count
$
        \beta_n:=\lfloor b_n^{1/d}\rfloor^d. 
$
Then, \(\beta_n\in\mathbb N\), \(\beta_n/b_n\to1\), and
\(r_n\beta_n^{1/d}\to0\).

The main result of this section is the LDP for the fixed radius component count.
\begin{theorem}  \label{t:LDP.Skn}
Under the radius \eqref{e:LDP.radius} satisfying \eqref{e:growth.an}, the sequence $(S_{k,n}/b_n)_{n\ge1}$ satisfies an LDP with speed sequence $(b_n)_{n\ge1}$ and  rate function 
\begin{equation}  \label{e:rate.func.Tkn}
I_k(x) := \begin{cases}
x\log (x/s_k) - x+s_k & \text{ if } x\ge0, \\
\infty & \text{ if } x<0, 
\end{cases}
\end{equation}
where
\[
s_k := \frac{\ta_d^{m_k}}{k!}\int_{(\R^d)^{k-1}} e^{-G(z_1,\dots,z_{k-1})} \dif \bz
\]
with 
$$
G(z_1,\dots, z_{k-1}) := \int_{\bigcup_{i=1}^{k-1}B_{\|z_i\|/2} (z_i/2)} \|x\|^{1-d}\dif x, \ \ \ z_i \in \R^d. 
$$
Equivalently, with \(\al_k\) as in Lemma~\ref{l:PY.pk}, we have
\(
        s_k={\al_k\ta_d^{m_k}}/{k}.
\)
\end{theorem}

\begin{remark}
The other  interesting setup is the LDP for 
$$
\widetilde S_{n,k} :=\sum_{\Y\subset \Pn, \, |\Y|=k} h(s_n^{-1}\Y) \one \big\{ (\Pn\setminus \Y) (\B(\Y; r_n))=0\big\}, 
$$
where $s_n:=r_n(nr_n^d)^{-1}\ll r_n$. We could apply the same approach as that for Theorem \ref{t:LDP.Skn}. 
\end{remark}

Before commencing the proof, we first observe that the speed sequence $(b_n)_{n\ge1}$ may exhibit either logarithmic or polynomial growth, depending on the choice of $a_n$. Specifically, if we take 
$$
a_n = \log n -(k-1)(d-1)\log \log n -\log \log \log n, 
$$
then, $b_n\sim \log \log n$ as $n\to\infty$. On the other hand, if we set $a_n=\log \log \log n$, then 
$$
b_n \sim \frac{n}{(\log \log n)(\log \log \log n)^{(k-1)(d-1)}}, \ \ \ n\to\infty. 
$$

For the proof of Theorem \ref{t:LDP.Skn}, we partition the flat torus $\bbT^d$ into $\beta_n$ sub-cubes $Q_1, \dots, Q_{\beta_n}$, with each of volume $\beta_n^{-1}$.  In particular, we set $Q_1=[0,\beta_n^{-1/d}]^d$. In addition, let us define 
$$
M_\ell(k)  := \big\{  x\in Q_\ell: \inf_{y\in \partial Q_\ell} \|x-y\| \le kr_n \big\}, 
$$
and consider 
$$
\eta_{k,n} = \sum_{\ell=1}^{\beta_n} \sum_{\Y\subset \Pnr_{\QmM(k)}} g_n(\Y, \Pn\setminus \Y) =: \sum_{\ell=1}^{\beta_n} \eta_{k,n}^{(\ell)}. 
$$
By the homogeneity of $\Pn$, the variables
$(\eta_{k,n}^{(\ell)})_{\ell=1}^{\beta_n}$ are identically distributed. Furthermore, they are
independent due to spatial independence of $\Pn$ and the following property: whenever $\Y\subseteq Q_\ell\setminus M_\ell(k)$, the
condition $g_n(\Y,\Pn\setminus\Y)=1$ is determined by the restriction of
$\Pn$ to $Q_\ell$.  

The first step identifies the limiting mean of a single block; this is the constant that later appears in the rate function. 
In what follows, $C$ as well as its variants represents a generic constant, which is independent of $n$. 

\begin{proposition}  \label{p:finite.mean}
For any choice of indices \(\ell=\ell(n)\in\{1,\ldots,\beta_n\}\),
\[
        \E[\eta_{k,n}^{(\ell)}]\to s_k,
        \qquad n\to\infty .
\]
\end{proposition}

The following result gives the high density Palm asymptotic needed to
identify the limiting block mean.  Let \(\mathcal C_1(0,\mathcal H_\lambda\cup\{0\})\) denote the vertex set of the component of
\(G(\mathcal H_\lambda\cup\{0\},1)\) containing the origin.  Define
\begin{equation}\label{e:def.pk.lambda}
p_k(\lambda):=\P\big(\big|\mathcal C_1(0,\mathcal H_\lambda\cup\{0\})\big|=k\big).
\end{equation}
For $r>0$, a locally finite point set $\xi$ in either $\bbT^d$ or
$\R^d$, and $x\in\xi$, let $\mathcal C_r(x,\xi)$ denote the vertex set
of the component of $G(\xi,r)$ containing $x$.
Then, by the multivariate Mecke formula for Poisson point processes (see, e.g., \cite[Chapter 4]{last:penrose:2017}), 
\begin{align}
p_k(\lambda) &= \E\Big[ \sum_{\Y\subset \mathcal H_\lambda, \, |\Y|=k-1} \one \big\{  G(\Y\cup \{0\}, 1) \text{ is a component of } G(\mathcal H_\lambda\cup \{0\}, 1)\big\}\Big]  \label{e:pk.lambda.Mecke} \\
&=\frac{\lambda^{k-1}}{(k-1)!} \int_{(\R^d)^{k-1}} h(0,y_1,\dots,y_{k-1}) \P \Big(\mathcal H_\lambda  \big( \B(\{0, y_1,\dots,y_{k-1}\};  1) \big)=0  \Big) \dif \by \notag \\
&= \frac{\lambda^{k-1}}{(k-1)!} \int_{(\R^d)^{k-1}} h(0,y_1,\dots,y_{k-1}) e^{-\lambda |\B (\{0,y_1,\dots,y_{k-1}\}; 1)|} \dif \by.\notag 
\end{align}

The asymptotic form of this Palm probability will be used throughout the paper.
\begin{lemma}[Theorem 2.1 in \cite{penrose:yang:2025}]  \label{l:PY.pk}
As $\lambda\to\infty$, 
$$
p_k(\lambda) = \al_k \lambda^{-m_k} e^{-\ta_d \lambda} \big( 1+O(\lambda^{-1}) \big), 
$$
where 
$$
\al_k:= \frac{1}{(k-1)!}\int_{(\R^d)^{k-1}}e^{-G(z_1,\dots,z_{k-1})}\dif \bz. 
$$
\end{lemma}

\begin{proof}[Proof of Proposition \ref{p:finite.mean}]
The proof reduces the block mean to the Palm probability in Lemma~\ref{l:PY.pk}; the only remaining issue is to show that the boundary effect of each  block is negligible.
By the multivariate Mecke formula, 
$$
\E[\eta_{k,n}^{(\ell)}] = \frac{n^k}{k!}\, \int_{(\QmM(k))^k} h_n(\bx)e^{-n|\B(\bx;  r_n)|}\dif \bx, 
$$
where $\bx=(x_1,\dots,x_k)\in (\R^d)^k$. 
Performing the change of variables $x_i=x+r_ny_{i-1}$, $i=1,\dots,k$ (with $y_0\equiv 0$), while using the location invariance of \eqref{e:def.hn} and the volume of  \eqref{e:def.union.ball}, 
\begin{align*}
\E[\eta_{k,n}^{(\ell)}] &=\frac{n(nr_n^d)^{k-1}}{k!}\int_{\QmM(k)}\int_{(\R^d)^{k-1}} h(0,\by) e^{-nr_n^d|\B(\{ 0,\by\}; 1)|} \prod_{i=1}^{k-1} \one\{ x+r_ny_i\in \QmM(k) \} \dif \by \dif x   \\
&= \frac{n(nr_n^d)^{k-1}}{k!} |\QmM(2k)| \int_{(\R^d)^{k-1}} h(0,\by) e^{-nr_n^d|\B(\{ 0,\by\}; 1)|} \dif \by \\
&\quad + \frac{n(nr_n^d)^{k-1}}{k!}\int_{(\QmM(k))\setminus( \QmM(2k))}\int_{(\R^d)^{k-1}} h(0,\by) e^{-nr_n^d|\B(\{ 0,\by\}; 1)|} \\
&\qquad \qquad \qquad \qquad \qquad\times \prod_{i=1}^{k-1} \one\{ x+r_ny_i\in \QmM(k) \} \dif \by \dif x \\
&=: A_n + B_n. 
\end{align*}
By Lemma \ref{l:PY.pk} and \eqref{e:pk.lambda.Mecke}, as $n\to\infty$, 
\begin{align*}
A_n = \frac{n}{k}  |\QmM(2k)| p_k(nr_n^d) \sim \frac{\al_k}{k} \ta_d^{m_k} n |\QmM(2k)|a_n^{-m_k}e^{-a_n}. 
\end{align*}
Since 
$$
b_nr_n^d
=
 n a_n^{-m_k}e^{-a_n}\frac{a_n}{n\ta_d}
=
\ta_d^{-1}a_n^{1-m_k}e^{-a_n}
\to0, \ \ \ n\to\infty, 
$$
we have $|\QmM(2k)| \sim |Q_\ell|=\beta_n^{-1}$ as $n\to\infty$, and since $b_n/\beta_n\to1$, it follows that $A_n\to s_k$.
It now remains to demonstrate that $B_n\to0$ as $n\to\infty$. Note first that 
$
B_n\le n \big| M_\ell(2k) \big| p_k (nr_n^d). 
$
Since $\big| M_\ell(2k) \big|  \le C r_n \beta_n^{-1+1/d}$, it can be  concluded from Lemma \ref{l:PY.pk} that 
$$
B_n \le C nr_n \beta_n^{-1+1/d} p_k(nr_n^d)  \le C (b_n/\beta_n)r_n\beta_n^{1/d} \to 0, \ \ \ n\to\infty. 
$$
\end{proof}

The next step shows that one block contains asymptotically a Poisson number of counted components.
\begin{proposition}  \label{p:total.variation.Poisson}
For every $\ell\ge1$, 
\begin{equation}  \label{e:total.variation.Poisson}
d_{\ms{TV}} \big( \mathcal L(\eta_{k,n}^{(\ell)}), \mathcal L(\ms{Poi} (s_k)) \big) := \sup_{A\in \B(\R)} \Big| \P(\eta_{k,n}^{(\ell)}\in A) - \P\big( \ms{Poi}(s_k) \in A \big) \Big| \to0,  \ \ \ n\to\infty, 
\end{equation}
where $\ms{Poi}(s_k)$ stands for a Poisson random variable with mean $s_k$. 
\end{proposition}

\begin{proof}
The proof uses a local-dependence Poisson approximation: after inserting a candidate component, only counted components whose exclusion regions overlap can affect the error.
Our proof is closely related to   \cite[Theorem 2.8]{iyer:yogeshwaran:2020} and \cite[Theorem 5.1]{penrose:yang:2025}, where the key technique employed  is Theorem 3.1 in \cite{penrose:2018}. 
For a $k$-point set $\Y\subset \R^d$ and a point measure $\xi$ in $\R^d$,  we set 
$$
f(\Y, \xi):=g_n(\Y, \xi) \one \{ \Y\subset \QmM(k) \}, 
$$
and 
$$
F(\xi) := \sum_{\Y\subset \xi, \, |\Y|=k} f(\Y, \xi\setminus \Y). 
$$
Then, $W_n:= F(\Pn) = \eta_{k,n}^{(\ell)}$. Moreover,  define for $\bx \in (\bbT^d)^k$,
$$
\Pn^\bx := \big\{  \Pn \cap \B(\bx;  r_n)^c\big\}\cup \bx, 
$$
and set 
$$
U_n^\bx := W_n, \ \ \text{ and } \ \ V_n^\bx := \sum_{\substack{\Y\subset \Pn^\bx, \\ \Y\neq \bx, \, |\Y|=k}} f(\Y, \Pn^\bx\setminus \Y). 
$$
For $\bx\in (\bbT^d)^k$ with $p_n(\bx):=\E[f(\bx, \Pn)]>0$, equivalently, $h_n(\bx)=1$ and $\bx\subset \QmM(k)$, 
we see that  
\begin{equation}  \label{e:1+V}
1+V_n^\bx=  f(\bx, \Pn^\bx\setminus \bx) + V_n^\bx= F(\Pn^\bx). 
\end{equation}
Moreover, if $f(\bx, \Pn)=1$, then $\Pn\big( \B(\bx; r_n) \big)=0$. Thus, $\Pn^\bx = \Pn \cup \bx$ and by \eqref{e:1+V}, 
$$
1+V_n^\bx \stackrel{d}{=} F(\Pn\cup \bx) \big|  \{ f(\bx, \Pn)=1\}. 
$$

Next, we decompose $U_n^\bx$ and $V_n^\bx$ as follows. 
\begin{align*}
U_n^\bx &= \sum_{\substack{\Y\subset \Pn, \\ \B(\bx;  r_n)\cap \B(\Y;  r_n)\neq \emptyset}} \hspace{-10pt}f(\Y, \Pn\setminus \Y) + \sum_{\substack{\Y\subset \Pn, \\ \B(\bx;  r_n)\cap \B(\Y;  r_n)= \emptyset}} \hspace{-10pt}f(\Y, \Pn\setminus \Y) =: W_n^{(1)}(\bx) + W_n^{(3)}(\bx), \\
V_n^\bx &= \sum_{\substack{\Y\subset \Pn^\bx, \, \Y \neq \bx \\ \B(\bx;  r_n)\cap \B(\Y; r_n)\neq \emptyset}} \hspace{-10pt}f(\Y, \Pn^\bx\setminus \Y) + \sum_{\substack{\Y\subset \Pn^\bx, \, \Y\neq \bx\\ \B(\bx; r_n)\cap \B(\Y; r_n)= \emptyset}} \hspace{-10pt}f(\Y, \Pn^\bx\setminus \Y) \\
&=: W_n^{(2)}(\bx) + \sum_{\substack{\Y\subset \Pn^\bx, \, \Y\neq \bx\\ \B(\bx;  r_n)\cap \B(\Y;  r_n)= \emptyset}} \hspace{-10pt}f(\Y, \Pn^\bx\setminus \Y).
\end{align*}
It is straightforward  to show  that 
$$
W_n^{(3)}(\bx) = \sum_{\substack{\Y\subset \Pn^\bx, \, \Y\neq \bx\\ \B(\bx;  r_n)\cap \B(\Y; r_n)= \emptyset}} \hspace{-10pt}f(\Y, \Pn^\bx\setminus \Y), 
$$
and hence, 
$$
|U_n^\bx-V_n^\bx| \le W_n^{(1)}(\bx) +  W_n^{(2)}(\bx). 
$$
Writing $w_n^{(i)}(\bx):=\E[W_n^{(i)}(\bx)]$ for $i=1,2$ and applying Theorem 3.1 in \cite{penrose:2018}, we have that 
$$
d_{\ms{TV}}\Big( \mathcal L(\eta_{k,n}^{(\ell)}), \mathcal L\big(\ms{Poi}(\E[\eta_{k,n}^{(\ell)}])\big) \Big) \le \frac{1\wedge \big( \E[\eta_{k,n}^{(\ell)}] \big)^{-1}}{k!}\, (I_{1,n}+I_{2,n}) \le I_{1,n}+I_{2,n}, 
$$
where, for $i=1,2$, 
$$
I_{i,n} = n^k \int_{(\bbT^d)^k} w_n^{(i)}(\bx) p_n(\bx) \dif \bx = n^k \int_{(\QmM(k))^k} w_n^{(i)}(\bx) h_n(\bx) e^{-n|\B (\bx; r_n)|}\dif \bx. 
$$
By Lemma 7.3 in \cite{bobrowski:adler:2014} and Proposition \ref{p:finite.mean}, 
$$
d_{\ms{TV}}\Big( \mathcal L\big(\ms{Poi}(\E[\eta_{k,n}^{(\ell)}])\big), \mathcal L\big(\ms{Poi}(s_k)\big) \Big) \le \big| \E[\eta_{k,n}^{(\ell)}] - s_k \big| \to 0, \ \ \text{as } n\to\infty. 
$$
Hence, \eqref{e:total.variation.Poisson} follows if one can show that $I_{1,n}+I_{2,n}\to0$ as $n\to\infty$. 

For the analysis of $I_{1,n}$, it follows from the multivariate Mecke formula and Proposition \ref{p:finite.mean} that 
\begin{align*}
I_{1,n} &=\frac{n^{2k}}{k!} \int_{(\QmM(k))^k} \int_{(\QmM(k))^k} h_n(\bx)h_n(\bz)\\
&\qquad \qquad \qquad\qquad \times \one \big\{ \B(\bx;  r_n)\cap \B(\bz; r_n)\neq \emptyset \big\} e^{-n \big\{|\B(\bx;  r_n)|+ |\B(\bz;  r_n)|\big\}}\dif \bx \dif \bz. 
\end{align*}
Since $| \B(\bx; r_n)| \ge \ta_d r_n^d$ and $| \B(\bz; r_n)| \ge \ta_d r_n^d$, one can see that as $n\to\infty$, 
\begin{align*}
I_{1,n} &\le e^{-2\ta_d nr_n^d} n^{2k} \int_{(\QmM(k))^k} \int_{(\QmM(k))^k} h_n(\bx) h_n(\bz)\one \big\{ \B(\bx;  r_n)\cap \B(\bz;  r_n)\neq \emptyset \big\} \dif \bx \dif \bz\\
&= e^{-2\ta_d nr_n^d} n (nr_n^d)^{2k-1}\big| \QmM(k) \big| \int_{(\R^d)^{2k-1}} h(0,y_1,\dots, y_{k-1})h(y_k,\dots, y_{2k-1})\\
&\qquad \qquad  \times \one \big\{ \B(\{0,y_1,\dots,y_{k-1}\}; 1)\cap \B(\{y_k,\dots,y_{2k-1}\}; 1) \big\} \prod_{i=1}^{2k-1} \one \{ x+r_ny_i\in \QmM(k) \}  \dif \by \\
&\le C e^{-2\ta_d nr_n^d} n (nr_n^d)^{2k-1}\beta_n^{-1 } \\
&=C (b_n/\beta_n)a_n^{m_k+2k-1} e^{-a_n} \to 0. 
\end{align*}

Subsequently, write for $\bx=(x_1,\dots, x_k) \in (\QmM(k))^k$, 
\begin{align} 
w_n^{(2)}(\bx) &= \sum_{j=0}^{k-1} \binom{k}{j} \E \bigg[ \hspace{-10pt} \sum_{\substack{\Z \subset \Pn \cap \B(\bx; r_n)^c, \, |\Z|= k-j, 
\\ \B(\{x_1,\dots,x_j,\Z\}; r_n) \cap \B(\bx; r_n)\neq \emptyset}} \hspace{-25pt}g_n\big( \{x_1,\dots,x_j, \Z\}, \Pn^\bx\setminus \{x_1,\dots,x_j, \Z\} \big)  \label{e:wn(2)} \\
&\qquad \qquad\qquad \qquad\qquad \qquad\qquad \qquad\times \one \{\Z \subset \QmM(k)\} \bigg]. \notag 
\end{align}
Then, 
\begin{align*}
&g_n\big( \{x_1,\dots,x_j, \Z\}, \Pn^\bx\setminus \{x_1,\dots,x_j, \Z\} \big) \\
&= h_n(x_1,\dots, x_j, \Z)\, \one \Big\{ (\Pn \setminus \Z) \big( \B(\{x_1,\dots,x_j,\Z\}; r_n)  \cap \B(\bx; r_n)^c\big)=0  \Big\}\\
&\qquad \qquad \qquad\qquad \qquad \times \one \big\{ x_i \in \B(\{x_1,\dots,x_j,\Z\}; r_n)^c, \, i=j+1,\dots,k  \big\}\\
&\le h_n(x_1,\dots, x_j, \Z)\, \one \Big\{ (\Pn \setminus \Z) \big( \B(\{x_1,\dots,x_j,\Z\}; r_n)  \cap \B(\bx; r_n)^c\big)=0  \Big\}. 
\end{align*}
Substituting this bound into $w_n^{(2)}(\bx)$ and applying the Mecke formula to \eqref{e:wn(2)}, we see that 
$$
I_{2,n}\le \sum_{j=0}^{k-1} \binom{k}{j}\frac{I_{2,j,n}}{(k-j)!}, 
$$
where 
\begin{align}
I_{2,j,n} &:=  n^{2k-j} \int_{\bx \in (\QmM(k))^k}\int_{(z_{j+1}, \dots, z_k)\in \big( (\QmM(k))\cap B(\bx,  r_n)^c \big)^{k-j}}\label{e:I2j} \\
&\qquad \qquad \qquad \qquad \qquad h_n(\bx) h_n(x_1,\dots,x_j,z_{j+1},\dots,z_{k}) \notag  \\
&\qquad \qquad \qquad \qquad \qquad\times \one \Big\{   \B( \bx; r_n) \cap \B \big(\{  x_1,\dots,x_j,z_{j+1},\dots,z_{k}\};  r_n \big)  \neq  \emptyset \Big\} \notag \\
&\qquad \qquad \qquad  \qquad \qquad\times e^{-n\big| \B(\{\bx,z_{j+1},\dots,z_k\}  ;r_n)  \big|} \dif z_{j+1} \cdots \dif z_{k}\dif \bx,  \notag 
\end{align}
where $\bx=(x_1,\dots,x_k)$. 

Suppose first that $j \in \{1, \dots, k-1\}$. In this case, the indicator function in \eqref{e:I2j} equals $1$.  
Furthermore, since $(z_{j+1}, \dots, z_k) \subset \B(\bx; r_n)^c$ implies $\|z_k - x_1\| > r_n$, we obtain
\begin{equation}  \label{e:lower.bdd.two.balls}
\big|\B(\{\bx,z_{j+1},\dots,z_k\}; r_n)  \big| \ge \big|B_{r_n}(x_1) \cup B_{r_n}(z_k) \big| \ge (1+\vep_0) \ta_d r_n^d
\end{equation}
for some $\vep_0>0$. Now, we have 
\begin{align*}
I_{2,j,n} &\le e^{-(1+\vep_0)\ta_d nr_n^d} n^{2k-j}  \int_{\bx \in (\QmM(k))^k}\int_{(z_{j+1}, \dots, z_k)\in ( (\QmM(k)) )^{k-j}} \hspace{-80pt} h_n(\bx)h_n(x_1,\dots,x_j,z_{j+1},\dots, z_k)  \dif z_{j+1} \cdots \dif z_{k}\dif \bx\\
&\le e^{-(1+\vep_0)\ta_d nr_n^d}  n (nr_n^d)^{2k-j-1}|Q_\ell| \int_{(\R^d)^{2k-j-1}}\hspace{-15pt} h(0,y_1,\dots, y_{k-1})h(0,y_1,\dots,y_{j-1}, y_k, \dots, y_{2k-j-1})\dif \by \\
&= C e^{-(1+\vep_0)\ta_d nr_n^d}   n (nr_n^d)^{2k-j-1}\beta_n^{-1} \\
&= C (b_n/\beta_n)a_n^{m_k+2k-j-1} e^{-\vep_0a_n} \to 0, \ \ \ n\to\infty. 
\end{align*}

Suppose next that $j=0$. Then, by using the bound \eqref{e:lower.bdd.two.balls} and repeating the same argument as above, we obtain that 
\begin{align*}
I_{2,0,n} &\le e^{-(1+\vep_0)\ta_d nr_n^d} n^{2k}  \int_{(\QmM(k))^k}\int_{( (\QmM(k)) )^{k}} h_n(\bx) h_n(\bz) \one \big\{\B(\bx; r_n)   \cap  \B(\bz; r_n) \neq \emptyset \big\}\dif \bz \dif \bx \\
&\le e^{-(1+\vep_0)\ta_d nr_n^d} n (nr_n^d)^{2k-1} |Q_\ell| \int_{(\R^d)^{2k-1}} h(0,y_1,\dots,y_{k-1}) h(y_k,\dots, y_{2k-1})\\
&\qquad \qquad \qquad \qquad \qquad\qquad  \times \one \big\{ \B(\{0,y_1,\dots,y_{k-1}\}; 1)  \cap \B(\{y_k,\dots,y_{2k-1}\}; 1)\neq \emptyset \big\}\dif \by \\
&= C (b_n/\beta_n)a_n^{m_k+2k-1} e^{-\vep_0 a_n}\to 0, \ \ \ n\to\infty. 
\end{align*}
\end{proof}

The same local counting estimate controls exponential moments within one block
and the components lost at the block boundaries.  We formulate it for a
deterministic set of possible component anchors, so that both applications
follow from the same coloring argument.  For a $k$-point set $\Y$, let
$\ell(\Y)$ be its lexicographically first point, and define
\[
        T_n(A):=\sum_{\substack{\Y\subset\Pn,\,|\Y|=k}}
                 g_n(\Y,\Pn\setminus\Y)\one\{\ell(\Y)\in A\}. 
\]
The following estimate is uniform over deterministic  sets. Given a Borel set $A$, let  
$$
A^{+u}:=\{x\in\bbT^d:\ms{dist}(x,A)\le u\}. 
$$
\begin{lemma}[Exponential moments for localized component counts]
\label{l:localized.exp.moments}
There are constants $C<\infty$ and $q=q(d,k)\in\mathbb N$ such that,
for all sufficiently large $n$, all Borel sets $A\subset\bbT^d$,
and all $t>0$,
\[
        \log\E [e^{tT_n(A)}]
        \le C(e^{qt}-1)b_n|A^{+Cr_n}|.
\]
\end{lemma}

\begin{proof}
Partition the torus into $N_n^d$ half-open cubes $J$ of side
$u_n=N_n^{-1}$, where $N_n=\lceil 4\sqrt d/r_n\rceil$.  For large $n$,
$cr_n\le u_n\le r_n/(4\sqrt d)$.  Define $Z_J:=T_n(J)$.
Since the diameter of $J$ is smaller than $r_n$, two distinct
components cannot have their anchors in the same $J$; thus,
$Z_J\in\{0,1\}$.  The Palm identity and Lemma~\ref{l:PY.pk} give
\[
        \P(Z_J=1)\le\E [Z_J]
        \le n|J|p_k(nr_n^d)\le Cb_nr_n^d.
\]
An internally connected $k$-point set with an anchor point in $J$ lies within
$(k-1)r_n$ of that anchor point, and its isolation can be decided by inspecting
points one radius further out.  Therefore, $Z_J$ depends only on
$\Pn\cap J^+$, where we have written $J^+:=J^{+(k+2)r_n}$ for simplicity.  These enlarged cubes have an
intersection graph of uniformly bounded degree.  Color the cubes with
$q=q(d,k)$ colors, giving classes
$\mathcal J_n^{(1)},\ldots,\mathcal J_n^{(q)}$, so that the variables
$Z_J$ within each color class are independent.

Let $\mathcal J_n(A)$ consist of the cubes meeting $A$.  These cubes lie
in $A^{+Cr_n}$ and have volume at least $cr_n^d$; hence,
\(
        \#\mathcal J_n(A)\le Cr_n^{-d}|A^{+Cr_n}|.
\)
Using $T_n(A)\le\sum_{J\in\mathcal J_n(A)}Z_J$, H\"older's inequality
across the color classes, and $\log(1+x)\le x$, we obtain
\[
\begin{aligned}
        \log\E [e^{tT_n(A)}]
        &\le\frac1q\sum_{i=1}^q
        \sum_{J\in\mathcal J_n^{(i)}\cap\mathcal J_n(A)}
        \log\big(1+(e^{qt}-1)\P(Z_J=1)\big)\\
        &\le C(e^{qt}-1)b_nr_n^d\#\mathcal J_n(A)
        \le C(e^{qt}-1)b_n|A^{+Cr_n}|.
\end{aligned}
\]
This proves the estimate.
\end{proof}

The localized estimate gives the exponential integrability required for the
blocked LDP.  It also permits passage from the distributional one-block
Poisson limit to convergence of exponential moments, since a slightly
larger exponent ensures uniform integrability.
\begin{lemma}[Uniform exponential moments for one block]
\label{l:one.block.exp.moments}
For every $a>0$,
\[
        \sup_{n\ge1}\E [e^{a\eta_{k,n}^{(1)}}]<\infty.
\]
Consequently, for every $a\in\R$,
\[
        \log\E [e^{a\eta_{k,n}^{(1)}}]\to  s_k(e^a-1).
\]
\end{lemma}

\begin{proof}
Every component counted by $\eta_{k,n}^{(1)}$ has its anchor point in $Q_1$,
so $\eta_{k,n}^{(1)}\le T_n(Q_1)$.  Since
$r_n\beta_n^{1/d}\to0$, we have
$|Q_1^{+Cr_n}|\le C\beta_n^{-1}$ for all large $n$.
Lemma~\ref{l:localized.exp.moments} therefore gives
\[
        \log\E [e^{a\eta_{k,n}^{(1)}}] 
        \le C(e^{qa}-1)b_n/\beta_n\le C_a.
\]
The finitely many remaining $n$ are harmless because the number of
components is bounded by the total of Poisson points.
For $a\le0$, convergence of the exponential moments follows directly
from Proposition~\ref{p:total.variation.Poisson}.  For $a>0$, apply the
bound above with $a'>a$ to obtain uniform integrability, and then use
the same Poisson approximation.  The limiting moment generating function
is $\exp\{s_k(e^a-1)\}$, as asserted.
\end{proof}

Once the one-block Laplace transform is controlled, the blocked approximation satisfies the desired LDP by a direct product argument.
\begin{proposition}
\label{p:LDP.eta}
The sequence \((\eta_{k,n}/b_n)_{n\ge1}\) satisfies an LDP with speed
\((b_n)_{n\ge1}\) and rate function \(I_k\) in \eqref{e:rate.func.Tkn}.
\end{proposition}

\begin{proof}
The proof is based on the G\"artner--Ellis theorem in large deviation theory: independence gives the limiting log-Laplace transform, and the endpoint \(x=0\) is checked separately.
By construction,
\(
        \eta_{k,n}=\sum_{\ell=1}^{\beta_n}\eta_{k,n}^{(\ell)},
\)
where the summands are i.i.d. Hence, for \(a\in\R\),
\[
        \Lambda_n(a)
        :=
        b_n^{-1}\log \E [e^{a\eta_{k,n}}]
        =
        \frac{\beta_n}{b_n}
        \log \E [e^{a\eta_{k,n}^{(1)}}].
\]
By Lemma~\ref{l:one.block.exp.moments}, we have
$\Lambda_n(a)\to \Lambda(a):=s_k(e^a-1)$ for $a\in\R$.
The Gärtner--Ellis theorem gives the LDP upper bound and the lower bound at
all \(x>0\), with
\[
        I_k(x)=\sup_{a\in\R}\{ax-s_k(e^a-1)\}.
\]
Thus, for \(x>0\), we have
\(
        I_k(x)=x\log(x/s_k)-x+s_k,
\)
 while \(I_k(x)=\infty\) for \(x<0\).

It remains to verify the lower bound at \(x=0\). If \(G\subset\R\) is open
and \(0\in G\), then, for all large \(n\),
\[
        \P(\eta_{k,n}/b_n\in G)
        \ge
        \P(\eta_{k,n}=0)
        =
        \P(\eta_{k,n}^{(1)}=0)^{\beta_n}.
\]
Since \(\eta_{k,n}^{(1)}\Rightarrow\ms{Poi}(s_k)\), we have
\(
        \P(\eta_{k,n}^{(1)}=0)\to e^{-s_k}.
\)
Together with \(\beta_n/b_n\to1\), this gives
\[
        \liminf_{n\to\infty}b_n^{-1}
        \log\P(\eta_{k,n}/b_n\in G)
        \ge
        -s_k=-I_k(0).
\]
This completes the LDP for \(\eta_{k,n}/b_n\).
\end{proof}

To complete the proof of Theorem~\ref{t:LDP.Skn}, we verify that
\((S_{k,n}/b_n)_{n\ge1}\) and \((\eta_{k,n}/b_n)_{n\ge1}\) are exponentially
equivalent. The blocking and exponential moment argument is adapted from the
dense regime nearest neighbor argument of \cite{hirsch:kang:owada:2023};
the one cube probability estimate below follows from the
high density component asymptotic of \cite{penrose:yang:2025}.
The following proposition establishes exponential equivalence.

\begin{proposition}
\label{p:exp.equiv.S.eta}
Let \(\delta>0\). Then,
\(
        \limsup_{n\to\infty} b_n^{-1}
        \log \P\big( |S_{k,n}-\eta_{k,n}|>\delta b_n \big)
        =
        -\infty .
\)
\end{proposition}

\begin{proof}
We first identify the components omitted by the blocking.  Let
\[
        M_n^*:=\bigcup_{i=1}^{\beta_n}
        \{x\in Q_i:\ms{dist}(x,\partial Q_i)\le2kr_n\}.
\]
If the anchor point $x=\ell(\Y)$ of a connected $k$-point set lies in
$Q_i\setminus M_n^*$, then every $y\in\Y$ lies within $(k-1)r_n$
of $x$.  Consequently,
\(
        \ms{dist}(y,\partial Q_i)
        >2kr_n-(k-1)r_n>kr_n.
\)
Thus, $\Y\subset Q_i\setminus M_i(k)$, and every component whose
anchor point lies outside $M_n^*$ is counted in $\eta_{k,n}$.  In particular,
\(
        0\le S_{k,n}-\eta_{k,n}\le T_n(M_n^*).
\)
The block geometry and $r_n\beta_n^{1/d}\to0$ give
\(
        |(M_n^*)^{+Cr_n}|\le Cr_n\beta_n^{1/d}=o(1).
\)
For each fixed $a>0$, Lemma~\ref{l:localized.exp.moments} now yields
\(
        \limsup_{n\to\infty}b_n^{-1}
        \log\E [e^{aT_n(M_n^*)}]\le0.
\)
By the exponential Markov inequality, for $\delta>0$,
\[
\begin{aligned}
        \limsup_{n\to\infty}b_n^{-1}
        \log\P(|S_{k,n}-\eta_{k,n}|>\delta b_n)
        \le-a\delta+
        \limsup_{n\to\infty}b_n^{-1}
        \log\E [e^{aT_n(M_n^*)}]
        \le-a\delta.
\end{aligned}
\]
Letting $a\to\infty$ proves exponential equivalence.
\end{proof}

\begin{proof}[Proof of Theorem~\ref{t:LDP.Skn}]
By Proposition~\ref{p:LDP.eta}, the sequence
\((\eta_{k,n}/b_n)_{n\ge1}\) satisfies an LDP with speed \(b_n\) and rate
function \(I_k\). By Proposition~\ref{p:exp.equiv.S.eta}, for every
\(\delta>0\),
\[
        \limsup_{n\to\infty} b_n^{-1}
        \log
        \P\big(
        |S_{k,n}-\eta_{k,n}|>\delta b_n
        \big)
        =
        -\infty .
\]
Thus, $(S_{k,n}/b_n)_{n\ge1}$ and $(\eta_{k,n}/b_n)_{n\ge1}$ are exponentially equivalent at
speed \(b_n\). The exponential equivalence theorem for large deviations
therefore implies that $(S_{k,n}/b_n)_{n\ge1}$ satisfies the same LDP as
$(\eta_{k,n}/b_n)_{n\ge1}$, namely with speed \(b_n\) and rate function \(I_k\) in
\eqref{e:rate.func.Tkn}.
\end{proof}

\bigskip

\section{Functional extremes for fixed size components}
\label{sec:extreme-k-clusters}

The preceding sections studied large deviations for the number of
\(k\)-components at a fixed threshold sequence.  A natural refinement is to
ask for a functional version, in whichdistinct thresholds are observed
simultaneously.  This is the analogue, for fixed size components, of the
functional large deviation viewpoint developed for \(k\)-nearest neighbour
volumes in \cite{hirsch:kang:owada:2023}.

For \(k\)-components, the corresponding scale variable is the external
isolation volume.  Roughly speaking, once a \(k\)-point set is internally connected,
its isolation volume records how far the surrounding ball can be expanded
before the first external point enters its neighbourhood.  Thus, instead of
only counting components at one prescribed threshold, we record the centered
thresholds at which base scale \(k\)-components cease to be isolated.

The main object of this section is the resulting
\(k\)-component isolation volume point process.  The tail measure theory
of Dombry, Tillier and Wintenberger
\cite{dombry:tillier:wintenberger:2022} provides the setting for our results.
We establish the required convergence directly through factorial moments.
A boundary layer estimate then transfers the point process limit to the
moving threshold component count path.

Throughout this section, \(\Pn\) is a homogeneous Poisson point process of
intensity \(n\) on the flat \(d\)-dimensional unit torus \(\bbT^d\).
Fix \(k\ge2\), recall \(m_k=(k-1)(d-1)\), and assume
\[
        a_n=\log n-m_k\log\log n+\delta_n,
        \qquad
        \delta_n\to\infty,
        \qquad
        \delta_n=o(\log n).
\]
The speed sequence in this setting  is defined by 
\[
        b_n=n a_n^{-m_k}e^{-a_n},
        \qquad
        r_n(t)=\Big(\frac{t}{n\ta_d}\Big)^{1/d},
        \qquad
        \rho_n=r_n(a_n).
\]
Then,
\(
        b_n=e^{-\delta_n}(1+o(1))\to0,
	\)
	and
	\(
        n\ta_d\rho_n^d=a_n.
\)
In contrast to Theorem~\ref{t:LDP.Skn}, where \(b_n\to\infty\), the present regime has \(b_n\to0\).
For a finite configuration \(\X\subset\bbT^d\) and a finite set
\(\Y\subset\X\), let
\[
\begin{aligned}
        r_{\rm conn}(\Y)&:=\inf\{r>0:\Y\text{ is connected in }G(\Y,r)\},\\
        r_{\rm ext}(\Y;\X)&:=\inf\{r>0:(\X\setminus\Y)(\B(\Y;r))>0\}.
\end{aligned}
\]
For a \(k\)-point set \(\Y\subset\X\), define the centered internal-connection level
and external isolation volume by 
\[
        B_n(\Y):=n\ta_d r_{\rm conn}(\Y)^d-a_n,
        \qquad
        D_n(\Y;\X):=n\ta_d r_{\rm ext}(\Y;\X)^d-a_n.
\]
Thus, \(B_n(\Y)\) is the centered level at which \(\Y\) becomes internally
connected, while \(D_n(\Y;\X)\) is the centered external isolation volume of
\(\Y\) relative to \(\X\).  The component count 
process is then defined by 
\[
        C_n^{(k)}(t)
        :=
        \sum_{\substack{\Y\subset\Pn\\ |\Y|=k}}
        \one\{B_n(\Y)\le t<D_n(\Y;\Pn)\},
        \qquad t\ge0.
\]
Equivalently, \(C_n^{(k)}(t)\) counts the \(k\)-point sets that are internally
connected at radius \(r_n(a_n+t)\) and have no external neighbour within the
same radius.  Moreover, let
\[
        E_n(\Y;\X):=
        \left\{
        r_{\rm conn}(\Y)\le\rho_n,
        \quad
        (\X\setminus\Y)(\B(\Y;\rho_n))=0
        \right\}.
\]
This is the event that \(\Y\) is connected at \(\rho_n\) and
has no external neighbour within that radius.  For \(h\in C_c((0,\infty))\), the space of continuous functions on $(0,\infty)$ with compact support, 
define the score
\[
        H_{n,h}(\Y;\X):=h(D_n(\Y;\X))\one_{E_n(\Y;\X)}.
\]
The associated \(k\)-component isolation-volume point process is
\[
        \Xi_n^{(k)}:=
        \sum_{\substack{\Y\subset\Pn\\ |\Y|=k}}
        \one_{E_n(\Y;\Pn)}\delta_{D_n(\Y;\Pn)}.
\]
On the exceptional event $|\Pn|=k$, the external radius is infinite.
We omit the corresponding atom at $+\infty$ when viewing
$\Xi_n^{(k)}$ as a point measure on $(0,\infty)$.  This event has
probability $e^{-n}n^k/k!=o(b_n^q)$ for every fixed $q$, and is accounted
for in the comparisons below.
The tail count $\Xi_n^{(k)}((t,\infty))$ records components already
internally connected at $\rho_n$ that remain externally isolated at
radius $r_n(a_n+t)$.  In contrast, $C_n^{(k)}(t)$ also counts components
that become internally connected between these two radii, provided
they remain externally isolated at the latter radius.  For test functions,
\(
        \Xi_n^{(k)}(h)=
        \sum_{\substack{\Y\subset\Pn\\ |\Y|=k}}H_{n,h}(\Y;\Pn).
\)

Let \(\mathcal N\) be the cone of point measures on \((0,\infty)\).  For
\(j\ge0\), let \(\mathcal N_j\subset\mathcal N\) be the cone of point
measures with at most \(j\) atoms.
Recall the Palm probability \(p_k(\lambda)\) from
\eqref{e:def.pk.lambda}.  We
use the Penrose--Yang asymptotic
\[
        p_k(\lambda)=\al_k\lambda^{-m_k}e^{-\ta_d\lambda}
        \bigl(1+O(\lambda^{-1})\bigr),
        \qquad \lambda\to\infty,
\]
and recall from Theorem~\ref{t:LDP.Skn} that
\(
        s_k={\al_k\ta_d^{m_k}}/{k}.
\)
For fixed $q\ge2$, define the measure 
\[
        \int F(\mu)\,d\nu_q^{(k)}(\mu)
        =
        \frac{s_k^q}{q!}
        \int_{(0,\infty)^q}
        F\Big(\sum_{i=1}^q\delta_{y_i}\Big)
        e^{-(y_1+\cdots+y_q)}\,dy_1\cdots dy_q
\]
for nonnegative measurable $F$.  More generally, for a nonnegative density $g$ on $(0,\infty)$ with finite
mass on every tail $(r,\infty)$, $r>0$, write
\[
        \int F(\mu)\,d\nu_q^g(\mu)
        =
        \frac1{q!}
        \int_{(0,\infty)^q}
        F\Big(\sum_{i=1}^q\delta_{y_i}\Big)
        \prod_{i=1}^q g(y_i)\,dy_i .
\]
Thus, $\nu_q^{(k)}=\nu_q^g$ for $g(y)=s_ke^{-y}$.
For the proof, we use the ordered factorial moment measures of the
isolation-volume process.  Put $N_n=\Xi_n^{(k)}((0,\infty))$ and, for
$m\ge1$, let $\alpha_{m,n}$ be the measure defined on Borel rectangles by
\[
        \alpha_{m,n}(A_1\times\cdots\times A_m)
        :=\E\bigg[\sum_{\Y_1,\ldots,\Y_m}^{\neq}
        \prod_{i=1}^m\one_{E_n(\Y_i;\Pn)}
        \one\{D_n(\Y_i;\Pn)\in A_i\}\bigg].
\]
Let $\lambda_k(dy):=s_ke^{-y}\,dy$ on $(0,\infty)$.
The next theorem establishes convergence at each fixed order of the
$k$-component isolation-volume point process.

\begin{theorem}
\label{t:DTW-k}
For every fixed $k\ge2$ and every fixed $q\ge2$,
\[
        b_n^{-q}\,\P(\Xi_n^{(k)}\in\cdot)
        \Longrightarrow
        \nu_q^{(k)}
        \qquad
        \text{in }\mathbb M(\mathcal N\setminus\mathcal N_{q-1}).
\]
\end{theorem}

The following proposition implies the theorem.
\begin{proposition}[Factorial moment convergence]
\label{p:k-factorial-measures}
For every fixed $m\ge1$,
\[
        b_n^{-m}\alpha_{m,n}
        \Longrightarrow \lambda_k^{\otimes m}
        \qquad\text{weakly as finite measures on }(0,\infty)^m.
\]
In particular,
\(
        b_n^{-m}\E\big[(N_n)_m\big]\to  s_k^m.
\)
\end{proposition}

The cases $m=q$ and $m=q+1$ of this proposition identify the
limiting $q$-atom configurations and show that configurations with more than
$q$ atoms are negligible at scale $b_n^q$.
\begin{proof}[Proof of Theorem~\ref{t:DTW-k}]
Fix $q\ge2$ and let
$F$ be bounded and continuous on
$\mathcal N\setminus\mathcal N_{q-1}$.  On $\{N_n=q\}$, each realization
of $\Xi_n^{(k)}$ is counted exactly $q!$ times by its ordered $q$-tuples.
On $\{N_n\ge q+1\}$, the error in replacing $F(\Xi_n^{(k)})$
by this ordered-tuple sum is bounded by
$C\|F\|_\infty(N_n)_{q+1}$.  Consequently,
\begin{align*}
&\Big|b_n^{-q}\E[F(\Xi_n^{(k)}) \one \{N_n\ge q\}]
  -\frac{b_n^{-q}}{q!}
    \int_{(0,\infty)^q}
       F\Big(\sum_{i=1}^q\delta_{y_i}\Big)\,
       \alpha_{q,n}(d\mathbf y)\Big| \\
  &\le C\|F\|_\infty b_n^{-q}\E\big[(N_n)_{q+1}\big]
   =O(b_n)\to 0,
\end{align*}
where the last estimate follows from
Proposition~\ref{p:k-factorial-measures} with $m=q+1$.
The map $\mathbf y\mapsto\sum_{i=1}^q\delta_{y_i}$ is continuous from
$(0,\infty)^q$ to $\mathcal N\setminus\mathcal N_{q-1}$.
By the same proposition with $m=q$, the integral tends to
\[
        \frac1{q!}\int_{(0,\infty)^q}
        F\Big(\sum_{i=1}^q\delta_{y_i}\Big)
        \prod_{i=1}^qs_ke^{-y_i}\,d\mathbf y
        =\nu_q^{(k)}(F).
\]
Thus, the normalized laws restricted to
$\mathcal N\setminus\mathcal N_{q-1}$ converge weakly as finite
measures to $\nu_q^{(k)}$.  This implies the asserted convergence in
$\mathbb M(\mathcal N\setminus\mathcal N_{q-1})$.
\end{proof}

\begin{remark}
The same argument also yields first order tail measure convergence
when \(q=1\). We state Theorem~\ref{t:DTW-k} for \(q\ge2\), where
interaction and higher-multiplicity estimates become essential.
\end{remark}

The point process limit also yields a functional statement for the
moving threshold component count.  We use the Skorokhod topology on
$\mathbb D_R:=D([0,R],\mathbb N_0)$, and set
\[
        \mathbb D_{R,j}:=\{f\in\mathbb D_R:\sup_{0\le t\le R}f(t)\le j\}.
\]
For $\mu\in\mathcal N$, define
\[
        \Phi_R(\mu)(t)=\mu((t,\infty)),\qquad 0\le t\le R,
\]
and put $\zeta_{q,R}^{(k)}:=\nu_q^{(k)}\circ\Phi_R^{-1}$.
The corollary follows from the theorem and a comparison of components
born at the base radius with those born afterward.

\begin{corollary}
\label{cor:k-functional-count}
For every fixed \(k\ge2\), every fixed \(q\ge2\), and every fixed
\(R<\infty\),
\[
        b_n^{-q}\,
        \P\bigl((C_n^{(k)}(t))_{0\le t\le R}\in\cdot\bigr)
        \Longrightarrow
        \zeta_{q,R}^{(k)}
        \qquad
        \text{in }\mathbb M(\mathbb D_R\setminus\mathbb D_{R,q-1}).
\]
\end{corollary}

We next control components that first become connected above the 
radius $\rho_n$.  For $R\ge0$, define
\[
        L_{n,R}
        :=\sum_{\substack{\Y\subset\Pn\\|\Y|=k}}
        \one\Big\{
        \rho_n<r_{\rm conn}(\Y)\le r_n(a_n+R),\,
        (\Pn\setminus\Y)(\B(\Y;\rho_n))=0
        \Big\}.
\]
The following estimate will be used in the moving threshold comparison
and in the proof of the one component asymptotics.

\begin{lemma}[Late-born components]
\label{lem:k-late-born}
For each fixed $0\le R<\infty$, there are constants $C<\infty$ and
$\varepsilon>0$ such that
\[
        \E [L_{n,R}]\le Cn a_n^{k-1}e^{-(1+\varepsilon)a_n}
\]
for all sufficiently large $n$.  In particular,
$\E [L_{n,R}]=o(b_n^q)$ for every fixed $q\ge1$.
\end{lemma}

Every late-born component contains an edge longer than $\rho_n$.
The resulting two-ball void bound gives the required exponential decay.
\begin{proof}
We estimate $\E [L_{n,R}]$ by the multivariate  Mecke formula.  Put
\(r_R:=r_n(a_n+R)\).  The insertion volume of a \(k\)-point set connected at radius
\(r_R\) is at most \(C r_R^{d(k-1)}\le C\rho_n^{d(k-1)}\).  On the additional
event \(r_{\rm conn}(\Y)>\rho_n\), the set \(\Y\) is connected in
\(G(\Y,r_R)\) but not in \(G(\Y,\rho_n)\).  Choose a spanning tree of
\(G(\Y,r_R)\).  This spanning tree cannot be contained entirely in
\(G(\Y,\rho_n)\), and therefore it contains an edge \(\{u,v\}\) with
\(
        \rho_n<|u-v|\le r_R .
\)
Since \(r_R/\rho_n\to1\), the scaled separation \(|u-v|/\rho_n\) lies in
\([1,2]\) for all large \(n\).  By continuity of
\(s\mapsto |B_1(0)\cup B_1(se_1)|\) and since this volume is strictly larger
than \(\ta_d\) on \([1,2]\), there is \(\varepsilon>0\) such that
\(
        |\B(\Y;\rho_n)|\ge (1+\varepsilon)\ta_d\rho_n^d
\)
for all large \(n\).  Therefore,
\begin{equation}\label{e:k-late-birth-estimate}
        \E [L_{n,R}]
        \le C n^k\rho_n^{d(k-1)}e^{-(1+\varepsilon)n\ta_d\rho_n^d}
        =C n a_n^{k-1}e^{-(1+\varepsilon)a_n}.
\end{equation}
Since \(b_n^q=n^q a_n^{-qm_k}e^{-q a_n}\),
\[
	b_n^{-q}\E [L_{n,R}]
        \le
        C n^{1-q}a_n^{k-1+qm_k}e^{(q-1-\varepsilon)a_n}
        =
        C a_n^{k-1+qm_k}
        e^{(q-1)(a_n-\log n)}e^{-\varepsilon a_n}
        \to0,
\]
because \(a_n-\log n=-m_k\log\log n+\delta_n=o(\log n)\) and
\(a_n\sim\log n\).
This proves the lemma.
\end{proof}

The late-birth estimate also controls the difference between the
actual moving threshold component count and its base born counterpart.
The next lemma makes this comparison uniform over a bounded time interval.
\begin{lemma}[Boundary-layer comparison for component counts]
\label{lem:k-count-comparison}
Let \(R<\infty\) and \(q\ge1\). Then,
\[
        b_n^{-q}
        \P\Big(
        \sup_{0\le t\le R}
        |C_n^{(k)}(t)-\Xi_n^{(k)}((t,\infty))|>0
	\Big)
        \to0 .
\]
\end{lemma}

The late-birth estimate leaves only the exceptional configuration to check.
\begin{proof}
For every fixed $q\ge1$, we have
\(
        \P(|\Pn|=k)=e^{-n}n^k/k!=o(b_n^q).
\)
Outside this event, we have $\Xi_n^{(k)}((t,\infty))\le C_n^{(k)}(t)$ for every $t\ge0$.
Any difference for $0\le t\le R$ requires a $k$-set counted by
$L_{n,R}$.  Hence, by Markov's inequality and
Lemma~\ref{lem:k-late-born},
\[
        \P\Big(\sup_{0\le t\le R}
          |C_n^{(k)}(t)-\Xi_n^{(k)}((t,\infty))|>0\Big)
        \le \E [L_{n,R}]+\P(|\Pn|=k)=o(b_n^q).
\]
\end{proof}

We now derive the functional corollary from the main theorem and the
comparison just proved.
\begin{proof}[Proof of Corollary \ref{cor:k-functional-count}]
Theorem~\ref{t:DTW-k} and the continuous-mapping theorem for tail
measures applied to \(\Phi_R\) give the asserted convergence for
\(\Phi_R(\Xi_n^{(k)})\).  The map \(\Phi_R\) is continuous
at \(\nu_q^{(k)}\)-almost every point measure, because under
\(\nu_q^{(k)}\) the atom locations have a joint density and hence there are
almost surely no atoms at \(0\), at \(R\), or at collision times.  The
comparison in Lemma~\ref{lem:k-count-comparison} transfers this
convergence to \(C_n^{(k)}\).
\end{proof}

\subsection{Proof of Proposition~\ref{p:k-factorial-measures}}

The factorial moment proposition is proved by expanding the ordered
factorial moments with the multivariate Mecke formula.  Well-separated
components contribute products of one component moments, whereas the
contributions of interacting components require an additional geometric
estimate.  We first establish the one component asymptotics using the
late-birth bound, then state the connected interaction estimate, and
finally derive the proposition from these two results.  The longer
geometric proof of the interaction estimate is given afterward.

The high density Palm asymptotic determines the one component tail.
\begin{lemma}
\label{lem:k-one-tail}
For every $R<\infty$ and every $h\in C_c((0,\infty))$, the following assertions hold.
\begin{enumerate}
\item[(i)] Uniformly for $0\le y\le R$,
\[
	\frac{\E\big[\Xi_n^{(k)}((y,\infty))\big]}{s_kb_ne^{-y}}\to1.
\]
\item[(ii)] One has
	$b_n^{-1}\E\big[\Xi_n^{(k)}(h)\big]\to s_k\int_0^\infty h(y)e^{-y}\,dy$.
\item[(iii)] There is $C<\infty$ such that, for all $T\ge0$,
	$\limsup_{n\to\infty}b_n^{-1}\E\big[\Xi_n^{(k)}((T,\infty))\big]\le Ce^{-T}$.
\end{enumerate}
\end{lemma}

The late-born component lemma makes the comparison with ordinary
$k$-components immediate.
\begin{proof}
Fix $R<\infty$ and $0\le y\le R$, and write $t=a_n+y$ and
$\lambda_t=t/\ta_d$.  The expected number of ordinary $k$-components
at radius $r_n(t)$ is $n p_k(\lambda_t)/k$.  By
Lemma~\ref{l:PY.pk}, uniformly in $y\in[0,R]$,
\[
        \frac{n}{k}p_k(\lambda_t)
        =s_kb_ne^{-y}(1+o(1)).
\]
Every component counted by $\Xi_n^{(k)}((y,\infty))$ is also counted
at radius $r_n(t)$.  Conversely, a component counted at that radius
but absent from $\Xi_n^{(k)}((y,\infty))$ is either born after
$\rho_n$ or belongs to the exceptional event $|\Pn|=k$.
The difference of expectations is therefore bounded by the expression
\(
        \E [L_{n,R}]+\P(|\Pn|=k)=o(b_n),
\)
by Lemma~\ref{lem:k-late-born} and
$\P(|\Pn|=k)=o(b_n)$.
This proves (i).

The normalized first-moment measures
$\mu_n=b_n^{-1}\E\big[\Xi_n^{(k)}(\cdot)\big]$ have tails converging to
$s_ke^{-y}$ for every $y\ge0$.  These tail limits give vague
convergence to $\lambda_k$, convergence of the total masses, and
therefore (ii).  Finally, for each $T\ge0$,
\[
        \E\big[\Xi_n^{(k)}((T,\infty))\big]
        \le\frac{n}{k}p_k((a_n+T)/\ta_d)
        \le C b_ne^{-T}
\]
for all sufficiently large $n$, by the high density Palm upper bound.
This proves (iii).
\end{proof}

The remaining estimate concerns connected groups of interacting
components.  It combines a spanning tree bound on the possible
locations with an additional void cost.  We state it here and prove
the factorial moment proposition before giving its geometric proof.
\begin{lemma}[Connected interaction estimate]
\label{lem:k-component}
Fix $L<\infty$ and $r\ge2$. Let $\mathcal A_C$ be the insertion domain corresponding to one connected component of size $r$ of the auxiliary interaction graph, with all test weights bounded by one. Then, for some $C<\infty$ and $\varepsilon>0$,
\[
        I_n(C):=
        \frac{n^{rk}}{(k!)^r}
        \int_{\mathcal A_C}
        \P\Big(\bigcap_{i=1}^rE_n(\Y_i(\mathbf x);\Pn\cup\{\mathbf x\})\Big)d\mathbf x
        \le C n^{rk}\rho_n^{d(rk-1)}e^{-(1+\varepsilon)a_n}=o(b_n^r).
\]
\end{lemma}

The one component asymptotics and the connected interaction estimate
are now in place.  We use them to control the separated and interacting
terms of the factorial moment expansion, with translation invariance
removing the near pair restriction from the product integrals.
\begin{proof}[Proof of Proposition~\ref{p:k-factorial-measures}]
The case $m=1$ follows from Lemma~\ref{lem:k-one-tail}, including
convergence of the total masses.  Fix $m\ge2$ and nonnegative
$h_1,\ldots,h_m\in C_c((0,\infty))$.  Choose $R<\infty$ such that their
supports lie in $(0,R)$, and fix $L>4$.  With $r_R=r_n(a_n+R)$, we have
$2r_R<L\rho_n$ for all sufficiently large $n$.

\medskip
\noindent\textbf{Step 1. Mecke expansion and interacting configurations.}
For $\mathbf x_i=(x_{i,1},\ldots,x_{i,k})$, let
$\Y_i=\{x_{i,1},\ldots,x_{i,k}\}$ and
$\mathbf x=(\mathbf x_1,\ldots,\mathbf x_m)$.
Two distinct $k$-sets satisfying $E_n$ cannot intersect: an internal
path from a shared vertex to a vertex outside the other set would give
that other set an external neighbour at distance at most $\rho_n$.
Thus, the ordered factorial sum runs over disjoint $k$-sets, and the
multivariate Mecke formula gives
\[
M_{m,n}(\mathbf h):=\int_{(0,\infty)^m}\prod_{i=1}^m h_i(y_i)\,\alpha_{m,n}(d\mathbf y)
=\frac{n^{km}}{(k!)^m}\int_{(\bbT^d)^{km}}
\E\bigg[\prod_{i=1}^m H_{n,h_i}\Big(\Y_i;\Pn\cup\bigcup_{j=1}^m\Y_j\Big)\bigg]d\mathbf x.
\]
Because $h_i$ vanishes outside $(0,R)$, its score is determined by the
Poisson process in the $r_R$-neighbourhood of $\Y_i$.  Join indices
$i$ and $j$ when $\ms{dist}(\Y_i,\Y_j)\le L\rho_n$.
Write $M_{m,n}^{\mathrm{int}}$ for the part of this integral where the
resulting interaction graph has an edge, and
$M_{m,n}^{\mathrm{sep}}$ for the edgeless part.

Let $\mathfrak G_m$ be the finite set of graphs on $\{1,\ldots,m\}$.
For $G\in\mathfrak G_m$, write $\mathcal B(G)$ for its connected
components.  For $B\in\mathcal B(G)$, let $J_{n,B}$ denote the
insertion integral involving only the clusters indexed by $B$, with
their interaction graph connected.  By Lemma~\ref{lem:k-one-tail},
$J_{n,\{i\}}=\E\big[\Xi_n^{(k)}(h_i)\big]=O(b_n)$; for $|B|=r\ge2$,
Lemma~\ref{lem:k-component} gives $J_{n,B}=o(b_n^r)$, with the bounded
test weights absorbed into the constant.  Clusters from distinct
components of $G$ are separated by more than $L\rho_n>2r_R$.
Their determining Poisson regions are therefore disjoint, so the
expectations factor.  Dropping the restrictions between these
components after factorization gives
\[
0\le M_{m,n}^{\mathrm{int}}(\mathbf h)
\le\sum_{\substack{G\in\mathfrak G_m\\E(G)\ne\varnothing}}
      \prod_{B\in\mathcal B(G)}J_{n,B}
=o(b_n^m).
\]
Indeed, every graph in the sum has a component of size at least two;
the sizes of its components add to $m$, and the number of graphs is
fixed.

\medskip
\noindent\textbf{Step 2. Separated configurations.}
Put $\Y(\mathbf y)=\{y_1,\ldots,y_k\}$ and
\[
\begin{aligned}
p_{n,i}(\mathbf y)
  &:=\E \Big[ H_{n,h_i}(\Y(\mathbf y);\Pn\cup\Y(\mathbf y))\Big],\\
A_{n,i}
  &:=\frac{n^k}{k!}\int_{(\bbT^d)^k}
        p_{n,i}(\mathbf y)\,d\mathbf y
   =\E\big[\Xi_n^{(k)}(h_i)\big].
\end{aligned}
\]
On the edgeless graph, the inserted sets and their $r_R$-neighbourhoods
are disjoint.  The corresponding expectation in the Mecke formula
therefore equals $\prod_{i=1}^m p_{n,i}(\mathbf x_i)$.
Each $p_{n,i}$ is translation invariant and vanishes unless its
$k$-set is connected at radius $\rho_n$.  If any $A_{n,i}$ vanishes,
both the separated integral and its product approximation vanish.
Otherwise, normalize $p_{n,i}(\mathbf y)d\mathbf y$ by its integral.
Under the product of these probability measures, the anchors
$x_{i,1}$ are independent and uniform on $\bbT^d$.

Every point of $\Y_i$ lies within $(k-1)\rho_n$ of $x_{i,1}$.
Thus, $\ms{dist}(\Y_i,\Y_j)\le L\rho_n$ implies that the
corresponding anchors are within $(L+2k-2)\rho_n$ of one another.
A union bound over $i<j$ shows that the complement of the edgeless
domain has probability at most $C\rho_n^d$ under the product measure.
Consequently,
\[
0\le \prod_{i=1}^m A_{n,i}-M_{m,n}^{\mathrm{sep}}(\mathbf h)
=\frac{n^{km}}{(k!)^m}\int_{\{\exists i<j:\,\ms{dist}(\Y_i,\Y_j)\le L\rho_n\}}
\prod_{i=1}^m p_{n,i}(\mathbf x_i)\,d\mathbf x
\le C\rho_n^d\prod_{i=1}^m A_{n,i}=o(b_n^m).
\]
The final estimate follows from Lemma~\ref{lem:k-one-tail}
and $\rho_n\to0$.

\medskip
\noindent\textbf{Step 3. Convergence of the factorial moment measures.}
The first two steps yield
$M_{m,n}(\mathbf h)=\prod_{i=1}^m A_{n,i}+o(b_n^m)$.
Lemma~\ref{lem:k-one-tail} therefore gives, for every choice of the
nonnegative test functions above,
\(
        b_n^{-m}M_{m,n}(\mathbf h)
        \to\prod_{i=1}^m\lambda_k(h_i).
\)
By linearity and approximation with finite sums of tensor products,
this proves vague convergence of $b_n^{-m}\alpha_{m,n}$ to
$\lambda_k^{\otimes m}$.

For the total masses, repeat Steps 1 and 2 with the unmarked
indicators $\one_{E_n(\Y_i;\Pn)}$.  They are determined by the
base-radius neighbourhoods, so the same factorization and bounds
apply.  The unmarked count can differ from $N_n$ only when
$|\Pn|=k$, because the atom at $+\infty$ is omitted from $\Xi_n^{(k)}$.
On this event, the unmarked count is at most one.  Thus, its
factorial moments of order $m\ge2$ agree with those of $N_n$, while
the one component expectations differ by at most
$\P(|\Pn|=k)=e^{-n}n^k/k!=o(b_n^m)$ for every fixed $m$.
Using $\E [N_n]\sim s_kb_n$ from Lemma~\ref{lem:k-one-tail}, we obtain
\[
        b_n^{-m}\E\big[(N_n)_m\big]
        =b_n^{-m}(\E [N_n])^m+o(1)
        \longrightarrow s_k^m
        =\lambda_k^{\otimes m}((0,\infty)^m).
\]
Vague convergence together with convergence of total masses implies
weak convergence as finite measures, as claimed.
\end{proof}

We conclude by proving the connected interaction estimate, using
the separation of
distinct components to obtain an additional exponential void cost.
\begin{proof}[Proof of Lemma \ref{lem:k-component}]
	We first establish the required geometric estimate. Namely, for every fixed $L<\infty$, there is $\varepsilon=\varepsilon(d,L)>0$ such that if two distinct counted $k$-point sets are at mutual distance at most $L\rho_n$, then
\(
        |\B(\Y;\rho_n)\cup\B(\Z;\rho_n)|
        \ge(1+\varepsilon)\ta_d\rho_n^d .
\)
Indeed, to see that this holds, let $E_n(\Y;\X)\cap E_n(\Z;\X)$ occur and assume
$\ms{dist}(\Y,\Z)\le L\rho_n$. Then, $\Y\cap\Z=\emptyset$.  The isolation part of $E_n(\Y;\X)$ gives
$(\X\setminus\Y)(\B(\Y;\rho_n))=0$.  Since
$\Z\subset\X\setminus\Y$, closed balls imply
$\ms{dist}(\Y,\Z)>\rho_n$.  Choose $u\in\Y$ and $v\in\Z$ with
$|u-v|=\ms{dist}(\Y,\Z)$.  Then, $\rho_n<|u-v|\le L\rho_n$.
Writing $s=|u-v|/\rho_n$, we have $s\in(1,L]$ and
\[
        |\B(\Y;\rho_n)\cup\B(\Z;\rho_n)|
        \ge
        |B_{\rho_n}(u)\cup B_{\rho_n}(v)|
        =
        \rho_n^d |B_1(0)\cup B_1(se_1)|.
\]
The continuous function $s\mapsto |B_1(0)\cup B_1(se_1)|$ has infimum
strictly larger than $\ta_d$ on $[1,L]$.  This gives the stated
$\varepsilon$.  Overlapping balls may have union volume smaller than
twice that of one ball; the argument uses only the uniform increase
over the volume of a single ball.

	Now we return to the proof of Lemma \ref{lem:k-component}. For inserted clusters
$\Y_1(\mathbf x),\ldots,\Y_r(\mathbf x)$, put an interaction edge between
$i$ and $j$ whenever we have
$\ms{dist}(\Y_i(\mathbf x),\Y_j(\mathbf x))\le L\rho_n$.  Thus, $\mathcal A_C$ is the set of inserted $rk$--tuples for which every
cluster is internally connected at radius $\rho_n$ and the auxiliary
interaction graph on the $r$ clusters is connected.  This is the domain
appearing in the definition of $I_n(C)$ in the statement.  The counted event
inside the probability also forces the inserted clusters to be pairwise
disjoint by the elementary no-overlap observation used above.
We bound the insertion volume first.  Choose a spanning tree of the
interaction graph.  The first anchor, say $x_{1,1}$, has torus volume $1$.
The remaining $k-1$ points of the first cluster are placed along an
internal spanning tree with edge lengths at most $\rho_n$, giving volume
at most $C\rho_n^{d(k-1)}$.  Each new cluster in the interaction spanning
tree has one anchor within distance $L\rho_n$ of a previously placed point,
giving volume at most $C\rho_n^d$, and its remaining $k-1$ points have
volume at most $C\rho_n^{d(k-1)}$.  Thus,
\[
        |\mathcal A_C|
        \le
        C\rho_n^{d(k-1)}
        \Big(\rho_n^d\rho_n^{d(k-1)}\Big)^{r-1}
        =
        C\rho_n^{d(rk-1)}.
\]
On the counted event, the Poisson process must avoid
$D_C(\mathbf x):=\bigcup_{i=1}^r\B(\Y_i(\mathbf x);\rho_n)$, up to
inserted points of Lebesgue measure zero.  Hence,
\(
        \P\Big(
        \bigcap_{i=1}^r E_n(\Y_i(\mathbf x);\Pn\cup\{\mathbf x\})
        \Big)
        \le
        e^{-n|D_C(\mathbf x)|}.
\)
Since the interaction graph is connected and $r\ge2$, some pair
$i,j$ is joined by an interaction edge.  The geometric estimate established at the
beginning of the proof gives
$|D_C(\mathbf x)|\ge(1+\varepsilon)\ta_d\rho_n^d$ on the counted event.
Therefore,
\[
        I_n(C)
        \le
        C n^{rk}\rho_n^{d(rk-1)}e^{-(1+\varepsilon)n\ta_d\rho_n^d}
        =
        C n^{rk}\rho_n^{d(rk-1)}e^{-(1+\varepsilon)a_n}.
\]
Because $\rho_n^d=a_n/(n\ta_d)$, we have
\(
        n^{rk}\rho_n^{d(rk-1)}=Cna_n^{rk-1}.
\)
Since $b_n^r=n^r a_n^{-rm_k}e^{-ra_n}$,
\[
        \frac{n^{rk}\rho_n^{d(rk-1)}e^{-(1+\varepsilon)a_n}}{b_n^r}
        \le
        C n^{1-r}a_n^{rk-1+rm_k}e^{(r-1-\varepsilon)a_n}.
\]
Now $a_n-\log n=-m_k\log\log n+\delta_n=o(\log n)$ and $\delta_n=o(\log n)$ and $a_n\sim\log n$. Thus, the last bound is
\[
        C a_n^{rk-1+rm_k}
        e^{(r-1)(-m_k\log\log n+\delta_n)}e^{-\varepsilon a_n}
        \to0,
\]
\end{proof}

%\bibliography{LDP-dense}

\begin{thebibliography}{10}

\bibitem{adams:konig:2008}
S.~Adams and W.~K\"onig.
\newblock Large deviations for many {B}rownian bridges with symmetrised
  initial-terminal condition.
\newblock {\em Probab. Theory Related Fields}, 142(1--2):79--124, 2008.

\bibitem{alexander:1993}
K.~S. Alexander.
\newblock Finite clusters in continuous percolation: compression and
  sphericality.
\newblock {\em Probab. Theory Related Fields}, 97:35--63, 1993.

\bibitem{andreis:konig:langhammer:patterson:2023}
L.~Andreis, W.~K\"onig, H.~Langhammer, and R.~I.~A. Patterson.
\newblock A large-deviations principle for all the components in a sparse
  inhomogeneous random graph.
\newblock {\em Probab. Theory Related Fields}, 186(1--2):521--620, 2023.

\bibitem{andreis:konig:patterson:2020}
L.~Andreis, W.~K\"onig, and R.~I.~A. Patterson.
\newblock A large-deviations principle for all the cluster sizes of a sparse
  {E}rd{\H{o}}s--{R}\'enyi graph.
\newblock {\em Random Structures Algorithms}, 59(4):522--553, 2021.

\bibitem{baccelli:blaszczyszyn:2009}
F.~Baccelli and B.~B{\l}aszczyszyn.
\newblock Stochastic geometry and wireless networks, volume {I}: Theory.
\newblock {\em Found. Trends Netw.}, 4:1--312, 2009.

\bibitem{bickel:breiman:1983}
P.~J. Bickel and L.~Breiman.
\newblock Sums of functions of nearest neighbor distances, moment bounds, limit
  theorems and a goodness of fit test.
\newblock {\em Ann. Probab.}, 11(1):185--214, 1983.

\bibitem{bobrowski:adler:2014}
O.~Bobrowski and R.~J. Adler.
\newblock Distance functions, critical points, and the topology of random
  \v{C}ech complexes.
\newblock {\em Homology Homotopy Appl.}, 16:311--344, 2014.

\bibitem{chatterjee:harel:2020}
S.~Chatterjee and M.~Harel.
\newblock Localization in random geometric graphs with too many edges.
\newblock {\em Ann. Probab.}, 48:574--621, 2020.

\bibitem{dombry:tillier:wintenberger:2022}
C.~Dombry, C.~Tillier, and O.~Wintenberger.
\newblock Hidden regular variation for point processes and the single/multiple
  large point heuristic.
\newblock {\em Ann. Appl. Probab.}, 32:191--234, 2022.

\bibitem{friedman:rafsky:1979}
J.~H. Friedman and L.~C. Rafsky.
\newblock Multivariate generalizations of the {W}ald--{W}olfowitz and {S}mirnov
  two-sample tests.
\newblock {\em Ann. Statist.}, 7(4):697--717, 1979.

\bibitem{hall:1986}
P.~Hall.
\newblock Clump counts in a mosaic.
\newblock {\em Ann. Probab.}, 14:424--458, 1986.

\bibitem{hall:1988}
P.~Hall.
\newblock {\em Introduction to the Theory of Coverage Processes}.
\newblock Wiley, New York, 1988.

\bibitem{henze:1988}
N.~Henze.
\newblock A multivariate two-sample test based on the number of nearest
  neighbor type coincidences.
\newblock {\em Ann. Statist.}, 16(2):772--783, 1988.

\bibitem{hirsch:kang:owada:2023}
C.~Hirsch, T.~Kang, and T.~Owada.
\newblock Large deviations for the volume of {$k$}-nearest neighbor balls.
\newblock {\em Electron. J. Probab.}, 28:1--27, 2023.

\bibitem{hirsch:otto:owada:thale:2025}
C.~Hirsch, M.~Otto, T.~Owada, and C.~Th\"ale.
\newblock Large deviations for the volume of hyperbolic {$k$}-nearest neighbor
  balls.
\newblock {\em Ann. Inst. Henri Poincar\'e Probab. Stat.}, 61(2):1443--1457,
  2025.

\bibitem{hirsch:owada:2023}
C.~Hirsch and T.~Owada.
\newblock Large deviation principle for geometric and topological functionals
  and associated point processes.
\newblock {\em Ann. Appl. Probab.}, 33:4008--4043, 2023.

\bibitem{iyer:yogeshwaran:2020}
S.~K. Iyer and D.~Yogeshwaran.
\newblock Thresholds for vanishing of `isolated' faces in random \v{C}ech and
  {V}ietoris--{R}ips complexes.
\newblock {\em Ann. Inst. Henri Poincar\'e Probab. Stat.}, 56:1869--1897, 2020.

\bibitem{jorritsma:komjathy:mitsche:2024}
J.~Jorritsma, J.~Komj\'athy, and D.~Mitsche.
\newblock Large deviations of the giant in supercritical kernel-based spatial
  random graphs.
\newblock {\em Probab. Theory Related Fields}, 2025.

\bibitem{last:penrose:2017}
G.~Last and M.~D. Penrose.
\newblock {\em Lectures on the Poisson Process}.
\newblock Cambridge University Press, 2017.

\bibitem{meester:roy:1996}
R.~Meester and R.~Roy.
\newblock {\em Continuum Percolation}.
\newblock Cambridge University Press, 1996.

\bibitem{penrose:2003}
M.~Penrose.
\newblock {\em Random Geometric Graphs, Oxford Studies in Probability 5}.
\newblock Oxford University Press, Oxford, 2003.

\bibitem{penrose:2018}
M.~Penrose.
\newblock Inhomogeneous random graphs, isolated vertices, and {P}oisson
  approximation.
\newblock {\em J. Appl. Probab.}, 55:112--136, 2018.

\bibitem{penrose:yang:2025}
M.~D. Penrose and X.~Yang.
\newblock On $k$-clusters of high-intensity random geometric graphs.
\newblock {arXiv:2209.14758}, 2025.

\end{thebibliography}
\end{document}